\documentclass{IEEEtran4PSCC}
\usepackage{algorithm, algpseudocode}
\usepackage{amsfonts}
\usepackage{xcolor}
\usepackage{caption, subcaption}

\ifCLASSINFOpdf
   \usepackage[pdftex]{graphicx}
\else
   \usepackage[dvips]{graphicx}
\fi
\usepackage[cmex10]{amsmath}
\makeatletter
\let\old@ps@headings\ps@headings
\let\old@ps@IEEEtitlepagestyle\ps@IEEEtitlepagestyle
\def\psccfooter#1{%
    \def\ps@headings{%
        \old@ps@headings%
        \def\@oddfoot{\strut\hfill#1\hfill\strut}%
        \def\@evenfoot{\strut\hfill#1\hfill\strut}%
    }%
    \def\ps@IEEEtitlepagestyle{%
        \old@ps@IEEEtitlepagestyle%
        \def\@oddfoot{\strut\hfill#1\hfill\strut}%
        \def\@evenfoot{\strut\hfill#1\hfill\strut}%
    }%
    \ps@headings%
}
\makeatother

\newcommand{\calT}{\mathcal{T}}
\newcommand{\calG}{\mathcal{G}}
\newcommand{\calB}{\mathcal{B}}

\newcommand{\calL}{\mathcal{L}}

\newcommand{\uon}{u^{on}}
\newcommand{\usu}{u^{su}}
\newcommand{\usd}{u^{sd}}

\newcommand{\uontg}{u^{on}_{t,g}}
\newcommand{\usutg}{u^{su}_{t,g}}
\newcommand{\usdtg}{u^{sd}_{t,g}}
\newcommand{\ptg}{p_{t,g}}
\newcommand{\qtg}{q_{t,g}}
\newcommand{\ptij}{p_{t,ij}}
\newcommand{\qtij}{q_{t,ij}}
\newcommand{\ptji}{p_{t,ji}}
\newcommand{\qtji}{q_{t,ji}}

\allowdisplaybreaks

\begin{document}
%
\title{On Solving Unit Commitment with Alternating Current Optimal Power Flow on GPU}

\author{\IEEEauthorblockN{Weiqi Zhang\IEEEauthorrefmark{1},
Youngdae Kim\IEEEauthorrefmark{2} and
Kibaek Kim\IEEEauthorrefmark{3}}
\IEEEauthorblockA{\IEEEauthorrefmark{1} Department of Chemical and Biological Engineering\\
University of Wisconsin-Madison, Madison, WI 53706}
\IEEEauthorblockA{\IEEEauthorrefmark{2} ExxonMobil Technology and Engineering Company -- Research, Annandale, NJ 08801}
\IEEEauthorblockA{\IEEEauthorrefmark{3} Mathematics and Computer Science Division\\
Argonne National Laboratory,
Lemont, IL 60439}
}

\maketitle

\begin{abstract} 
We consider the unit commitment (UC) problem that employs the alternating current optimal power flow (ACOPF) constraints, which is formulated as a mixed-integer nonlinear programming problem and thus challenging to solve in practice. We develop a new scalable algorithm based on alternating direction method of multiplier (ADMM), which enables quickly finding a good feasible solution of the UC-ACOPF problem. Our algorithm employs the component-based decomposition that solves a large number of independent small subproblems, each of which represents an individual grid component (e.g., generator, bus, and transmission line), as well as the UC subproblem. We implement the algorithm in Julia, as part of the existing package \texttt{ExaAdmm.jl}, which can efficiently run on GPUs as well as CPUs. The numerical results are reported by using IEEE test instances.
\end{abstract}


\thanksto{\noindent This material is based upon work supported by the U.S. Department of Energy, Office of Science, under contract number DE-AC02-06CH11357.
This research was conducted when the first two authors were affiliated with Argonne National Laboratory.}

\section{Introduction}
The unit commitment (UC) problem is concerned with finding a cost-minimizing production and operating schedule of generators of a power system to satisfy the demand over a time horizon, typically with hourly intervals, while taking into account various operational characteristics of generators such as ramp limits, start-up/shutdown times, and minimum on/off limits.
The UC problem is in general solved with approximate linear network constraints instead of more accurate nonconvex nonlinear alternating current (AC) network constraints, resulting in a mixed-integer linear programming (MILP).
This is mainly because its combinatorial nature combined with the AC network constraints will lead to a mixed-integer nonlinear programming (MINLP), which is computationally even more challenging to solve.

However, as pointed in~\cite{sauer2014} and demonstrated in~\cite{castillo2016}, incorporating AC network constraints into the UC problem can potentially provide a significant merit in terms of feasibility and its economic value.
Along this line of direction, some existing works in the literature  try to reflect the AC network constraints by employing relaxations, such as direct-current (DC) OPF relaxation \cite{cain2012history}, second-order conic relaxation \cite{kocuk2016strong, constante2021ac}, and semi-definite programming relaxation \cite{bai2009semi}.
But, it is still desired to directly solve the UC-ACOPF without relaxation, since it is often not trivial to recover AC feasible solutions from solutions of relaxation models \cite{liu2018global} and it can lead to a more economical schedule~\cite{castillo2016}.

In this paper we directly embed AC network constraints in the UC problem and introduce a decomposition scheme based on alternating direction method of multipliers (ADMM) to efficiently solve the resulting MINLP problem, called the UC-ACOPF.
Our scheme enables us to decompose the UC-ACOPF into the UC and ACOPF subproblems so that we can separately deal with mixed-integer and nonconvex nonlinear AC constraints.
Each of the UC and ACOPF subproblems can be further decomposed into grid components by time, resulting in a large number of independent small subproblems amenable to parallel computing.
Since the computation time of these subproblems is critical to the overall computational performance of our method, we exploit the massive parallel computing capability of GPUs for a significant acceleration of its computation time compared with that on CPUs.

In particular we exploit GPUs for both UC and ACOPF subproblems.
Our decomposition scheme allows us to move all the continuous variables into the ACOPF subproblems, thus making the UC subproblems consist purely of independent binary quadratic programs defined for each generator.
Instead of using the off-the-shelf solver to solve these problems, we show that a dynamic programming (DP) algorithm  can fully utilize its problem structure by efficiently performing a series of many array operations on GPUs.
In the case of the ACOPF subproblems, we present our formulations that leverage the existing batch solvers, ExaTron.jl~\cite{kim2021leveraging} and ExaAdmm.jl~\cite{kim2022accelerated}, in order to significantly accelerate the computation time of many small nonconvex nonlinear continuous optimization problems on GPUs.

The rest of the paper is organized as follows.
In Section~\ref{sec:uc-acopf_formulation} we present our UC-ACOPF formulation.
We describe three major elements of it:  the UC problem, the ACOPF problem, and the part that couples the UC problem with the ACOPF problem.
Section~\ref{sec:decomposition} presents our ADMM-based decomposition scheme together with formulations of subproblems and an in-depth description of our dynamic programming algorithm for the UC subproblems.
Numerical results are presented in Section~\ref{sec:exp} over IEEE test instances on Nvidia GPUs.
Section~\ref{sec:conclusion} concludes this paper with some discussion about future work.

\section{UC-ACOPF Formulation}
\label{sec:uc-acopf_formulation}

This section presents our UC-ACOPF model following the formulation in \cite{tejada2019unit} and \cite{constante2021ac}. We first give a brief overview of the UC-ACOPF model and then present its decomposable structure for the algorithm development.

The variables of the UC-ACOPF model can be partitioned into two parts, $x_t^\text{OPF}$ and $x_t^\text{UC}$, which contain values for defining ACOPF and UC problems for each time period $t \in \calT$, respectively:
\begin{align*}
    x_t^\text{OPF} &:= ([\ptg, \qtg]_{g\in\calG}, [\ptij, \qtij]_{(i,j)\in\calL},\\
    &\quad\quad [w^R_{t,i}]_{i\in\calB}, [w^R_{t,ij}, w^I_{t,ij}]_{(i,j)\in\calL}) \\
    x_t^\text{UC} & := ([\uontg, \usutg, \usdtg]_{g\in\calG}),
\end{align*}
where $p_{t,g}$ and $q_{t,g}$ are variables for real and reactive power generated at generator $g\in\calG$, $p_{t,ij}$ and $q_{t,ij}$ are real and reactive power flows over line $(i,j) \in \calL$, $w^R_{t,i},w^R_{t,ij},w^I_{t,ij}$ are for power flow equations between line $(i,j) \in \calL$ using a rectangular form of voltages at bus $i\in\calB$, and $\uontg, \usutg$, and $\usdtg$ represent binary variables for UC for generator $g\in\calG$.
We assume that $x_0^\text{OPF}$ and $x_0^\text{UC}$ are given, for example, from the previous system operations.
Then UC-ACOOF can be written as the following generic formulation:
\begin{subequations}
\label{eq:uc-acopf}
    \begin{align}
        \min_{x^\text{OPF},x^\text{UC}} \quad & \sum_{t\in\calT} f_t^\text{OPF}(x_t^\text{OPF}) + f^\text{UC}(x^\text{UC})\\
        \text{s.t.} \quad 
        & x^\text{OPF}_t \in \text{ACOPF}_t \quad \forall t\in\calT, \label{eq:uc-acopf-acopf}\\
        & x^\text{UC} := (x^\text{UC}_1,\dots,x^\text{UC}_T) \in \text{UC}, \label{eq:uc-acopf-uc}\\
        & R_{g,t}(x^\text{OPF}_{g,t},x^\text{OPF}_{g,t-1},x^\text{UC}_{g,t},x^\text{UC}_{g,t-1}) \leq 0 \quad \forall g\in\calG, t\in\calT,
    \end{align}
\end{subequations}
where $f_t^\text{OPF}$ and $f^\text{UC}$ are the corresponding cost functions, $\text{ACOPF}_t$ represents the set of feasible ACOPF solutions at time period $t\in\calT$,  UC represents the set of feasible UC solutions, and $R_t(\cdot)$ represents the generator bound constraints and the ramp-limit constraints that couple consecutive time periods for each $t\in\calT$.

\subsection{ACOPF Formulation}

We present the detailed formulation of ACOPF by defining the cost function $f_t^\text{OPF}$ and the feasible region ACOPF$_t$ given in \eqref{eq:uc-acopf-acopf}. We first define the objective cost function as follows:
\begin{align*}
    f_t^\text{OPF}(x_t^\text{OPF}) := \sum_{g\in\calG} f_{t,g}^\text{OPF}(p_{t,g}),
\end{align*}
which represents the sum of the total generation cost at each time $t\in\calT$.
The set ACOPF$_t$ of feasible OPF solutions is defined by the following set of constraints representing power flow equations, line limits, and bounds on voltage and generator variables:
\begin{subequations}
    \begin{align}
        & \sum_{g \in \calG_i} \ptg - P_{t,i} = G^s_i w^R_{t,i} + \sum_{j \in \calB_i}\ptij \quad \forall i \in \calB, \label{uc_acopf_p_balance} \\
        & \sum_{g \in \calG_i} \qtg - Q_{t,i} = -B^s_i w^R_{t,i} + \sum_{j \in \calB_i}\qtij \quad \forall i \in \calB, \label{uc_acopf_q_balance} \\
        & (\ptij)^2+(\qtij)^2 \leq \Bar{r}_{ij}^2 \quad \forall (i,j) \in \calL, \label{uc_acopf_ij_limit} \\
        & (\ptji)^2+(\qtji)^2 \leq \Bar{r}_{ij}^2 \quad \forall (i,j) \in \calL, \label{uc_acopf_ji_limit} \\
        & \ptij = G_{ii}w_{t,i}^R + G_{ij}w_{t,ij}^R + B_{ij}w_{t,ij}^I \quad \forall (i,j) \in \calL, \\
        & \qtij = -B_{ii}w_{t,i}^R - B_{ij}w_{t,ij}^R + G_{ij}w_{t,ij}^I \quad \forall (i,j) \in \calL, \\
        & \ptij = G_{jj}w_{t,j}^R + G_{ji}w_{t,ij}^R - B_{ji}w_{t,ij}^I \quad \forall (i,j) \in \calL, \\
        & \qtij = -B_{jj}w_{t,j}^R - B_{ji}w_{t,ij}^R - G_{ji}w_{t,ij}^I \quad \forall (i,j) \in \calL, \\
        & (w_{t,ij}^R)^2 + (w_{t,ij}^I)^2 = w_{t,i}^Rw_{t,j}^R \quad \forall (i,j) \in \calL, \\
        & \theta_{t,i} - \theta_{t,j} = \arctan (w_{t,ij}^I/w_{t,ij}^R) \quad \forall (i,j) \in \calL, \\
        & \underline{V}_i^2 \leq w_{t,i}^R \leq \Bar{V}_i^2 \quad \forall i \in \calB, \label{uc_acopf_w_bound} \\
        & 0 \leq \ptg \leq \Bar{P}_g \quad \forall g \in \calG, \label{acopf_gen_p_bounds} \\
        & 0 \leq \qtg \leq \Bar{Q}_g \quad \forall g \in \calG. \label{acopf_gen_q_bounds}
    \end{align}
    \label{eq:acopf}
\end{subequations}
Note that the constraints are independent of the UC variables, which are coupled to the ACOPF problem through the constraints described in \ref{s:coupling}.

\subsection{UC Formulation}
\label{s:uc-formulation}
Similar to the ACOPF, the objective function for the UC part is a sum of individual cost for each time period and generator:
\begin{align*}
    f^\text{UC}(x^\text{UC}) := \sum_{t\in\calT} \sum_{g\in\calG} f^\text{UC}_{t,g}(x^\text{UC}_g),
\end{align*}
where $x^\text{UC}_g := (\uontg, \usutg, \usdtg)$, $f^\text{UC}_{t,g}$ consists of the operation cost, start-up cost, and shutdown cost.
The constraint set UC in \eqref{eq:uc-acopf-uc} follows the standard UC formulation defined for each generator $g\in\calG$, which we denote by UC$_g$. For each $g\in\calG$, UC$_g$ is defined by the following constraints:
\begin{subequations}
\label{set:UC}
\begin{align}
& \sum_{t=1}^{L_g} (1-\uon_{t,g}) = 0, \label{uc_initial_on}\\
& \sum_{t=1}^{F_g} \uon_{t,g} = 0, \label{uc_initial_off}\\
& \uon_{t-1,g} - \uon_{t,g} + \usu_{t,g} - \usd_{t,g} = 0 \quad \forall t \in \calT\backslash\{1\}, \label{uc_acopf_gen_state} \\
& \sum_{i=t-T_g^U+1}^{t} \usu_{i,g} \leq \uon_{t,g} \quad \forall t \in \{T_g^U, ..., T\}, \label{uc_minup}\\
& \sum_{i=t-T_g^D+1}^{t} \usd_{i,g} \leq 1-\uon_{t,g} \quad \forall t \in \{T_g^D, ..., T\}, \label{uc_mindn}\\
& \uon_{t,g}, \usu_{t,g}, \usd_{t,g} \in \{0,1\} \quad \forall t \in \calT. \label{uc_acopf_binary}
\end{align}
\end{subequations}
The generator's initial status is given by constraints \eqref{uc_initial_on} and \eqref{uc_initial_off}.
Constraint \eqref{uc_acopf_gen_state} models the logic of turning on and off the generator.
If a generator is turned on (or off), then it has to stay on (or off) at least for a given number of time periods. These are modeled by constraints \eqref{uc_minup} and \eqref{uc_mindn}.
Note that with the definition of the objective function $f^\text{UC}$, the UC problem is separable for each generator $g\in\calG$.

\subsection{Coupling Constraints}
\label{s:coupling}
Here we present the coupling constraints $R_{g,t}(\cdot)$ that represent the following constraints:GAIL - this seems awkward -- why not just say we present the coupling constraints $R...$ aand then present them?
\begin{subequations}
    \begin{align}
        & \underline{P}_g \uontg \leq \ptg \leq \Bar{P}_g \uontg \label{uc_acopf_gen_p_bounds} \\
        & \underline{Q}_g \uontg \leq \qtg \leq \Bar{Q}_g \uontg \label{uc_acopf_gen_q_bounds} \\
        & \ptg - p_{t-1,g} \leq R_g^U \uon_{t-1,g} + S_g^U \usutg, \label{uc_acopf_gen_ramp_up_ineq} \\
        & \ptg - p_{t-1,g} \geq -R_g^D \uontg - S_g^D \usdtg.\label{uc_acopf_gen_ramp_down_ineq}
    \end{align}
    \label{eq:coupling}
\end{subequations}
Note that \eqref{uc_acopf_gen_p_bounds}--\eqref{uc_acopf_gen_q_bounds} couples ACOPF with UC.
Constraints \eqref{uc_acopf_gen_ramp_up_ineq}--\eqref{uc_acopf_gen_ramp_down_ineq} decide whether the ramping is in a normal running mode or startup/shutdown mode based on the unit commitment variables values.

\section{Decomposition Method for UC-ACOPF}
\label{sec:decomposition}

In this section we present the decomposition of UC-ACOPF into the UC and multiperiod ACOPF subproblems.
The multiperiod ACOPF subproblem is further decomposed into single-period ACOPF subproblems via temporal decomposition, each of which can be subsequently solved by the component-based decomposition method in \cite{mhanna2018adaptive,kim2021leveraging,kim2022accelerated}.

To apply our decomposition scheme, we reformulate the problem into a form amenable to the ADMM by introducing duplicate variable $\bar{x}=(\bar{x}^\text{UC},\bar{x}^\text{OPF})$ and extending $x^\text{OPF}$ to incorporate slack variables for the equality reformulation of~\eqref{eq:coupling}:
\begin{subequations}
    \label{eq:uc-acopf-ext}
    \begin{align}
        \min_{x^\text{OPF},x^\text{UC},\bar{x}} \quad & \sum_{t\in\calT} f_t^\text{OPF}(x_t^\text{OPF}) + f^\text{UC}(x^\text{UC})\\
        \text{s.t.} \quad 
        & \eqref{eq:acopf}-\eqref{set:UC} \notag\\
        & A_1x^\text{UC} + B_1\bar{x}^\text{UC}=0, \label{eq:uc_uc_coupling}\\
        & A_2x^\text{OPF} + B_2\bar{x}^\text{UC}=0, \label{eq:acopf_uc_coupling}\\
        & A_3x^\text{OPF} + B_3\bar{x}^{\text{OPF}}=0, \label{eq:acopf_acopf_coupling}
    \end{align}
\end{subequations}
where~\eqref{eq:uc_uc_coupling} represents the coupling between $x^\text{UC}$ and its duplicate $\bar{x}^\text{UC}$, \eqref{eq:acopf_uc_coupling} contains the coupling between $x^\text{OPF}$, and $\bar{x}^\text{UC}$ via the equality equivalent of~\eqref{eq:coupling}, and~\eqref{eq:acopf_acopf_coupling} encapsulates the coupling between the ACOPF components (generators, lines, and buses) as in~\cite{mhanna2018adaptive, kim2022accelerated}.
Specifically,~\eqref{eq:uc_uc_coupling}--\eqref{eq:acopf_uc_coupling} take the following form:
\begin{subequations}
    \begin{align}
        x^\text{UC}_{t,g} - \bar{x}^{\text{UC}}_{t,g} = 0, \label{eq:uc_duplicate}\\
        p_{t,g} - s^l_{t,p_g} - \underline{P}_g \bar{u}^\text{on}_{t,g} =0, \label{eq:pg_lo}\\
        p_{t,g} + s^u_{t,p_g} - \overline{P}_g\bar{u}^\text{on}_{t,g} =0,\\
        q_{t,g} - s^l_{t,q_g} - \underline{Q}_g\bar{u}^\text{on}_{t,g}=0,\\
        q_{t,g} + s^u_{t,q_g} - \overline{Q}_g\bar{u}^\text{on}_{t,g}=0,\\
        p_{t,g}-p_{t-1,g} - s^l_{t,r_g} + R^D_g\bar{u}^{\text{on}}_{t-1,g}+S^D_g\bar{u}^{\text{su}}_{t,g}=0,\\
        p_{t,g}-p_{t-1,g} + s^u_{t,r_g} - R^U_g\bar{u}^{\text{on}}_{t-1,g}-S^U_g\bar{u}^{\text{su}}_{t,g}=0. \label{eq:ramp_up}
    \end{align}    
\end{subequations}
In~\eqref{eq:pg_lo}-[\eqref{eq:ramp_up}, a slack variable $s$ has been introduced in order to equivalently reformulate~\eqref{eq:coupling} into equality constraints.
For~\eqref{eq:acopf_acopf_coupling}, we refer to~\cite{mhanna2018adaptive, kim2022accelerated}.

\subsection{Two-Level ADMM Method}

We apply the two-level ADMM method~\cite{sun2021two} to solve problem~\eqref{eq:uc-acopf-ext}, which provides theoretical convergence guarantees of ADMM iterates for nonconvex problems.
The basic idea of the two-level ADMM is to introduce an artificial variable $z$ with constraint $z=0$ to each of the coupling constraints like~\eqref{eq:uc_uc_coupling}-\eqref{eq:acopf_acopf_coupling} and to make it play a role of the last block of the coupling constraints so that it can resolve conflicting conditions~\cite{sun2021two} that prevent convergence of the ADMM algorithm.
The value of variable $z$ is guided to move toward 0 by applying an augmented Lagrangian to constraint $z=0$ via so-called outer-level iterations, and each outer-level iteration involves the solve of the problem containing coupling constraints with $z$ via ADMM.
Each iteration of the ADMM is called an inner-level iteration.

Following the steps in~\cite{sun2021two}, we add an artificial variable $z$ to the coupling constraints with constraint $z=0$ and apply the augmented Lagrangian relaxation with respect to constraint $z=0$, which results in the following:
\begin{subequations}
    \label{eq:outer}
    \begin{align}
        \min_{x^\text{OPF},x^\text{UC},\bar{x},z} \
        & \sum_{t\in\calT} f_t^\text{OPF}(x_t^\text{OPF}) + f^\text{UC}(x^\text{UC}) + \lambda^T z + \frac{\beta}{2}\|z\|^2 \notag\\
        \text{s.t.} \
        & \eqref{eq:acopf}-\eqref{set:UC} \notag\\        
        & A_1x^\text{UC} + B_1\bar{x}^\text{UC}+z_1=0, \label{eq:coupling_with_z1}\\
        & A_2x^\text{OPF} + B_2\bar{x}^\text{UC}+z_2=0, \label{eq:coupling_with_z2}\\
        & A_3x^\text{OPF} + B_3\bar{x}^\text{OPF}+z_3=0, \label{eq:coupling_with_z3}
    \end{align}
\end{subequations}
This problem serves as the outer-level iteration of the two-level ADMM method; at each outer-level iteration we solve~\eqref{eq:outer} and update $\lambda$ and $\beta$ by following the augmented Lagrangian algorithm until convergence.

For a given $(\lambda, \beta)$ we solve~\eqref{eq:outer} using ADMM, and this solve is called inner-level iterations.
To apply ADMM, an augmented Lagrangian with respect to~\eqref{eq:coupling_with_z1}-\eqref{eq:coupling_with_z3} is formed with constraints~\eqref{eq:acopf}--\eqref{set:UC}:
\begin{align*}
    & L_{\beta,\rho}(x^\text{OPF},x^\text{UC},\bar{x}^\text{OPF},\bar{x}^\text{UC},z,y,\lambda) \\
    &:= \sum_{t\in\calT} f_t^\text{OPF}(x_t^\text{OPF}) + f^\text{UC}(x^\text{UC}) + \lambda^T z + \frac{\beta}{2}\|z\|^2 \\
    & + y_1^T(A_1x^\text{UC}+B_1\bar{x}^\text{UC}+z_1) + \frac{\rho_1}{2}\|A_1x^\text{UC}+B_1\bar{x}^\text{UC}+z_1\|^2\\
    & + y_2^T(A_2x^\text{OPF}+B_2\bar{x}^\text{UC}+z_2) + \frac{\rho_2}{2}\|A_2x^\text{OPF}+B_2\bar{x}^\text{UC}+z_2\|^2\\
    & + y_3^T(A_3x^\text{OPF}+B_1\bar{x}^\text{OPF}+z_3) + \frac{\rho_3}{2}\|A_3x^\text{OPF}+B_1\bar{x}^\text{OPF}+z_3\|^2.
\end{align*}
Each ADMM iteration alternatively updates $x=(x^\text{OPF},x^\text{UC})$ and $\bar{x}=(\bar{x}^\text{OPF},\bar{x}^\text{UC})$ while keeping other variables fixed.
Specifically, at the $l$th inner-level ADMM iteration we compute the following:
\begin{subequations}
\begin{align}
    (x^\text{UC})^{l+1} &\gets \underset{x^\text{UC}\in\eqref{set:UC}}{\text{argmin}} \,L_{\beta^k,\rho}((x^\text{OPF})^l,x^\text{UC},(\bar{x})^l, z^l,y^l,\lambda^k) \label{step:uc}\\
    (x^\text{OPF})^{l+1} &\gets \underset{x^\text{OPF}\in\eqref{eq:acopf}}{\text{argmin}}\,L_{\beta^k,\rho}(x^\text{OPF},(x^\text{UC})^{l},\bar{x}^l, z^l,y^l,\lambda^k) \label{step:opf}\\
    (\bar{x}^\text{UC})^{l+1} &\gets \underset{\bar{x}^\text{UC}\in [0,1]^{\text{card}(\bar{x}^\text{UC})}}{\text{argmin}} \,L_{\beta^k,\rho}((x)^{l+1},\bar{x}^\text{UC}, z^l,y^l,\lambda^k) \label{step:uc_bar}\\
    (\bar{x}^\text{OPF})^{l+1} &\gets \underset{\bar{x}^\text{OPF}\in \eqref{eq:acopf}}{\text{argmin}} \,L_{\beta^k,\rho}((x)^{l+1},(\bar{x}^\text{UC})^l,\bar{x}^\text{OPF}, z^l,y^l,\lambda^k) \label{step:opf_bar}\\
    z^{l+1} &\gets \underset{z}{\text{argmin}}\, L_{\beta^k,\rho}((x)^{\ell+1},(\bar{x})^{l+1},z,y^l,\lambda^k) \label{step:z}\\
    y^{l+1} &\gets y^l + \rho (A (x)^{l+1} + B (\bar{x})^{l+1} + z^{l+1}) \label{step:y}.
\end{align}
\end{subequations}
We note that for a given $\bar{x}$ the UC and ACOPF problems can be completely separated.
Therefore,~\eqref{step:uc} and~\eqref{step:opf} can be solved in parallel.
The same applies to~\eqref{step:uc_bar} and~\eqref{step:opf_bar} as well.
Also, when we solve~\eqref{step:opf} and~\eqref{step:opf_bar}, the component-based decomposition scheme~\cite{kim2022accelerated} is applied.

Once we solve~\eqref{eq:outer}, 
the multiplier $\lambda^k$ and parameter $\beta^k$ are updated as follows:
\begin{align}
    \label{step:outer}
    \lambda^{k+1} &\gets \Pi_{[\underline{\lambda},\overline{\lambda}]}(\lambda^k + \beta^k z^{k})\\
    \beta^{k+1} &\gets \left\{
        \begin{array}{ll}
            \tau\beta^k &\text{if}\, \|z\|^k > \theta\|z\|^{k-1}\\
            \beta^k &\text{otherwise}
        \end{array}\right.,
\end{align}
where $z^k$ is the final value of $z$ when the inner level terminates, and we set $\tau=6$ and $\theta=0.8$.
The two-level ADMM method can be summarized as follows:
\begin{algorithm}
    \caption{Two-level ADMM for UC-ACOPF}
    \label{alg:basic}
    \begin{algorithmic}[1]
        \Require{UC-ACOPF data, $\rho$, $\epsilon$}
        \State Initialize $\lambda^0,\beta^0$ and $k\gets 0$.
        \While{$\|z\|>\epsilon$}
            \State Initialize $(x^\text{UC})^0,(\bar{x}^\text{OPF})^0,z^0,y^0$ and $l\gets 0$.
            \While{inner iteration not converged}
                \State Solve \eqref{step:uc} for $x^\text{UC}$. \label{step:uc}
                \State Solve \eqref{step:opf} for $x^\text{OPF}$. \label{step:opf}
                \State Solve \eqref{step:uc_bar} for $\bar{x}^\text{UC}$. \label{step:baruc}
                \State Solve \eqref{step:opf_bar} for $\bar{x}^\text{OPF}$. \label{step:baropf}
                \State Solve \eqref{step:z} for $z$. \label{step:z}
                \State Update multiplier $y$ by \eqref{step:y} and $l \gets l+1$.
            \EndWhile
            \State Update multiplier $\lambda$ and penalty $\beta$ by \eqref{step:outer}.
            \State Update $k \gets k+1$.
        \EndWhile
    \end{algorithmic}
\end{algorithm}

Algorithm \ref{alg:basic} assumes that the problem data, penalty parameter $\rho$, and tolerance $\epsilon$ are given.
The algorithm starts with initializing the values of $\lambda$ and $\beta$ and terminates when $\|z\| \leq \epsilon$ at the outer level.
Each outer-level iteration initializes the inner-level iterates $(x^\text{UC})^l,(\bar{x}^\text{OPF})^l,z^l,y^l$. For outer-level iteration $k>0$, the inner-level iterates are initialized from the previous outer-level iteration in order to exploit the warm start.
The inner level terminates if some termination criteria are met (e.g., $\epsilon$-stationary point of the outer-level problem \eqref{eq:outer} \cite{sun2021two}).
Because lines \ref{step:opf} and \ref{step:baropf} of Algorithm~\ref{alg:basic} are the same as in \cite{kim2022accelerated} and because line \ref{step:baruc} is straightforward, we focus on how line \ref{step:uc} can be efficiently solved via reformulation into dynamic programming and GPUs.

\subsection{Dynamic Programming for the UC Subproblem}

In Section \ref{s:uc-formulation} we present the set of constraints for the UC subproblem in \eqref{set:UC}, which is defined for each generator $g\in\calG$, such that
\begin{align*}
    & L_{\beta^k,\rho}((\bar{x}^\text{OPF})^l,x^\text{UC},(\bar{x})^l,z^l,y^l,\lambda^k) \\
    & \quad = \sum_{g\in\calG} L_{g,\beta^k,\rho}((\bar{x}_g^\text{OPF})^l,x_g^\text{UC},(\bar{x}_g)^l,z_g^l,y_g^l,\lambda^k).
\end{align*}
Moreover, at each inner-level iteration $l$, the objective function of the UC subproblem is also separable for each generator $g\in\calG$.
As a result, the UC subproblem for generator $g\in\calG$ can be written as
\begin{align}
    \min_{x_g^\text{UC}\in\text{UC}_g} \ L_{g,\beta^k,\rho}((\bar{x}_g^\text{OPF})^l,x_g^\text{UC},(\bar{x}_g)^l, z_g^l,y_g^l,\lambda^k). \label{uc_subproblem}
\end{align}
This problem is a binary quadratic programming problem that can be solved by an off-the-shelf solver on the CPU. However, we notice that problem \eqref{uc_subproblem} has a special structure that can be exploited for dynamic programming more efficiently than the existing solvers. For our purpose, there are two benefits of doing so. First, such specialized solution methods are able to utilize the structure and thus, as we will show later, are able to achieve a computational complexity much faster than a general IP solver. Second, a dynamic programming method can be implemented in the GPU as it is simply array operations, which eliminates the need for calling IP solvers on the CPU, and the associated communication between the CPU and GPU. 

We first restate problem \eqref{uc_subproblem} as follows: \textit{For each generator $g$, given an initial condition (on/off status $\uon_{0,g}$), minimum on and off time requirements $T^U_g$, $T^D_g$, and a series of cost functions $\{f^g_t\}_{t\in\mathcal{T}}$, what is the optimal trajectory $\{\uon_{t,g}\}_{t\in\mathcal{T}}$ over the time horizon $\mathcal{T}$}? Note that without loss of generality, we ignore the initial on/off time time requirement $F_g$/$L_g$ captured in constraints \eqref{uc_initial_on} and \eqref{uc_initial_off}. Moreover, we ignore the decision variables $\usu, \usd$ because their values can be inferred from the values of $\uon$: 
$$\usu_{t,g} = \max\{0, \uon_{t,g} - \uon_{t-1,g}\}$$
$$\usd_{t,g} = \max\{0, \uon_{t-1,g} - \uon_{t,g}\}$$.
As a result, the objective function of \eqref{uc_subproblem}  depends only on the current and previous states.
To simplify the description of dynamic programming algorithm, we denote such objective function by $L^\text{UC}_{g,t}(\uon_{t-1,g}, \uon_{t,g})$ for the rest of this section.

We describe our dynamic programming approach for solving each UC subproblem \eqref{uc_subproblem}. 
The main idea is to break down the problem into smaller simple decisions whether to stay or change the current status with respect to the cost function $L^\text{UC}_{g,t}(\uon_{t-1,g}, \uon_{t,g})$ by storing their optimal solution trajectory (i.e., memoization). 
Let $s\in\{0,1\}$ be the state of the generator representing on and off, and let $c_{t,s}$ and $u_{t,s,t:T}$ be the optimal cost and the optimal solution trajectory, respectively, of the subproblem starting at the beginning of time period $t$ with the previous state $s$. 

At the beginning of time period $t$, the generator needs to decide whether to stay at the current state $s$ or change to the other one $1-s$.
If the decision is made to stay at the current state $s$, the DP epoch moves to the next time period with the cost value
\begin{subequations}
    \label{keep_formula}
    \begin{equation}
        \label{c_keep_formula}
        c^\text{stay}_{t,s} := L^\text{UC}_{g,t}(s, s) + c_{t+1,s}
    \end{equation}
    with the corresponding solution trajectory to be
    \begin{equation}
        \label{u_keep_formula}
        u^\text{stay}_{t,s,t'} := 
        \begin{cases}
        s, &\quad t' = t \\
        u_{t+1,s,t'} &\quad t' \in [t+1, T].\\
        \end{cases}
    \end{equation}
\end{subequations}

If the decision is made to switch to the other state $1-s$, then we need to consider the minimum on/off time requirement based on the type of switching (i.e., on vs. off). Let $T^\text{min}_g$ be the minimum number of time periods required for the generator $g$ to stay at the new state. For example, $T^\text{min}_g=T^U_g$ if the new state is on; $T^\text{min}_g=T^D_g$ otherwise. The cost value of making this choice is 
\begin{subequations}
    \label{switch_formula}
    \begin{equation} 
    \label{c_switch_formula}
    \begin{split}
        c^\text{switch}_{t,s} & := L^\text{UC}_{g,t}(s, 1-s)
        + \sum_{t'=t+1}^{t+T^\text{min}_g-1} L^\text{UC}_{g,t'}(1-s, 1-s) \\ 
        & \quad + c_{t+T^\text{min}_g,1-s}
    \end{split}
    \end{equation}
    with the corresponding solution trajectory to be
    \begin{equation}
    \label{u_switch_formula}
        u^\text{switch}_{t,s,t'} :=
        \begin{cases}
        1-s, &\quad t' \in [t, t+T^\text{min}-1] \\
        u_{t+T^\text{min},1-s,t'} &\quad t' \in [t+T^\text{min}, T].\\
        \end{cases}
    \end{equation}
\end{subequations}
The optimal decision at time period $t$ is made by comparing the two cost values $c^\text{stay}_{t,s}$ and $c^\text{switch}_{t,s}$. For example, if $c^\text{stay}_{t,s} \leq c^\text{switch}_{t,s}$, then the decision is to stay with the optimal cost and trajectory of $c^\text{stay}_{t,s}$ and $u^\text{stay}_{t,s,t:T}$. The values of $c_{t,s}$ and $u_{t,s,t:T}$ will be updated with these results.

As the decision at each time period depends on the optimal solution of the following time periods, we can leverage a backward induction implementation of the dynamic programming method. Specifically, we start from the last time period $T$ and compute the optimal solution at each time period $t$ backward to the first time period. The process ends after making the decision for the first time period. Based on this idea, we propose the dynamic programming algorithm for solving UC subproblems, as described in Algorithm \ref{dp_algo}.

\begin{algorithm}[!ht]
\caption{DP for UC Subproblem of Generator $g$}\label{dp_algo}
\begin{algorithmic}[1]
    \Require{UC data, cost function $L^\text{UC}_{g,t}$}
    \For{$ t \in \{T, T-1, ..., 1\}$}
        \If{$\uon_{0,g} == 1$ and $t \leq T^U_g$}
            \State $u_{1,\uon_{0,g},1:T} \leftarrow u_{L_g+1, \uon_{0,g}, :}$
            
            \State $c_{1,\uon_{0,g}} \leftarrow \sum_{t=1}^{L_g} L^\text{UC}_{g,t}(\uon_{0,g}, \uon_{0,g}) + c_{L_g+1,\uon_{0,g}}$

            \State \textbf{break}
        \EndIf
        \If{$\uon_{0,g} == 0$ and $t \leq T^D_g$}
            \State $u_{1,\uon_{0,g},1:T} \leftarrow u_{F_g+1, \uon_{0,g}, :}$
            
            \State $c_{1,\uon_{0,g}} \leftarrow \sum_{t=1}^{F_g} L^\text{UC}_{g,t}(\uon_{0,g}, \uon_{0,g}) + c_{F_g+1,\uon_{0,g}}$ 

            \State \textbf{break}
        \EndIf

        \For{$s \in \{0, 1\}$}
            \State Compute $c^\text{stay}_{t,s}$ and $u^\text{stay}_{t,s}$ by \eqref{keep_formula}

            \State Compute $c^\text{switch}_{t,s}$ and $u^\text{switch}_{t,s}$ by \eqref{switch_formula}

            \If{$c^\text{stay}_{t,s} \leq c^\text{switch}_{t,s}$}
                \State Set $c_{t,s} \gets c^\text{stay}_{t,s}$
                and $u_{t,s, t:T} \gets u^\text{stay}_{t,s,t:T}$.
            \Else
                \State Set $c_{t,s} \gets c^\text{switch}_{t,s}$
                and $u_{t,s, t:T} \gets u^\text{switch}_{t,s,t:T}$.
            \EndIf
        \EndFor
    \EndFor
    \State \textbf{return} Optimal schedule $u_{1,u^{on}_{0,g}, 1:T}$ with the cost $c_{1, \uon_{0,g}}$
\end{algorithmic}
\end{algorithm}
The time complexity of this algorithm is $O(T)$, since it simply performs simple additions over $T$ time periods. The space complexity of this algorithm is $O(T^2)$, as it stores an optimal trajectory over time periods $t:T$ for each possible start time $t$. This shows that our algorithm solves the UC subproblem in linear time, theoretically way more efficiently than general-purpose integer programming solvers.

\section{Numerical Experiments}
\label{sec:exp}

\begin{table*}[!t]
    \centering
    \begin{tabular}{l|rrrrrrrrrrr}
             & \multicolumn{3}{c}{Parameters} &           & \multicolumn{2}{c}{Iterations} & Primal & \multicolumn{2}{c}{Time (s)} & \\
        Case & $\rho_\text{pq}$ & $\rho_\text{va}$ & $\rho_\text{uc}$ & Objective & Outer level & Total & Infeasibility & CPU & GPU & Speedup \\
        \hline
        case9   & 5e+3 & 1e+4 & 1e+4 & 9.1291e+4 & 9  & 320  & 1.8e-3 & 4.1   & 9.9 & 0.4 \\ 
        case30  & 5e+5 & 1e+6 & 1e+6 & 1.0603e+4 & 13 & 386  & 3.8e-3 & 18.3  & 18.9 & 0.9 \\
        case118 & 5e+4 & 1e+5 & 1e+5 & 2.1730e+6 & 12 & 456  & 8.2e-3 & 107.0 & 34.6 & 3.1 \\
        case300 & 5e+3 & 1e+4 & 1e+4 & 1.2210e+7 & 40 & 1961 & 1.2e-2 & 652.8 & 140.1 & 4.7 
    \end{tabular}
    \caption{Numerical results from the component-based ADMM algorithm for solving the test instances.}
    \label{tab:admm}
\end{table*}

We demonstrate the computational performance of our algorithm with the numerical experiments on IEEE test instances. 

\subsection{Experiment Setting}

The component-based ADMM algorithm proposed in this paper is implemented in an open-source Julia package \texttt{ExaAdmm.jl}. This package implements the component-based ADMM for solving single-period ACOPF, the numerical results of which has been demonstrated in \cite{kim2022accelerated}.
While leveraging the existing structure and design of the package particularly for the ACOPF subproblem, our algorithm implements the modified generator subproblems and the DP algorithm (Algorithm \ref{dp_algo}) for the UC subproblem.
Moreover, the DP algorithm is also written for both the CPU and GPU.
As in \cite{kim2021leveraging}, we solve the line subproblems by using the \texttt{ExaTron.jl} package for the GPU solver customized for a bound-constrained nonlinear optimization with up to 32 variables only.

The termination conditions of the algorithm are the same as in \cite{sun2021two,kim2022accelerated}. Moreover, similar to \cite{kim2022accelerated}, we use different values of the penalty parameter $\rho$ for the UC variables, the real/reactive power variables, and the other variables, denoted by $\rho_\text{uc}$, $\rho_\text{pq}$, and $\rho_\text{va}$, respectively. The specific values are reported later for each instance.
The ACOPF variables are initialized for the cold start, as in \cite{kim2022accelerated}. Specifically, the real and reactive power generation and voltage magnitude are initialized to the medium of their lower and upper bounds. Voltage angles are initialized to zero while fixing the reference voltage angle to zero. The power flow values are computed by using the initial voltage values.
The UC variables are initialized from a warm start procedure, where the multiperiod ACOPF model is solved, and UC variables are assigned 0 or 1 depending on the generation levels from the warm start solution.
All experiments were run on a workstation equipped with Nvidia’s Quadro GV100 and Intel Xeon 6140 CPU@2.30GHz.

\subsection{Test Instances}

We use 9-, 30-, 118-, and 300-bus IEEE test instances from \texttt{MATPOWER}. Based on the demand level provided in the original \texttt{.mat} files, we generated a demand trajectory for 24 hours from a time series of scaling factors and then multiplied by a discount factor of 0.7 to induce optimal solution with some generators being off (otherwise the demand level is too high, and we always see optimal solution with all generators on). The scaling factors were obtained from ISO-NE historical demand data in October 2018. The ramp of each generator was set to 10\% of its generation capacity. 

\subsection{Performance of Dynamic Programming}

\begin{figure}[!ht]
    \centering
    \includegraphics[width=0.48\textwidth]{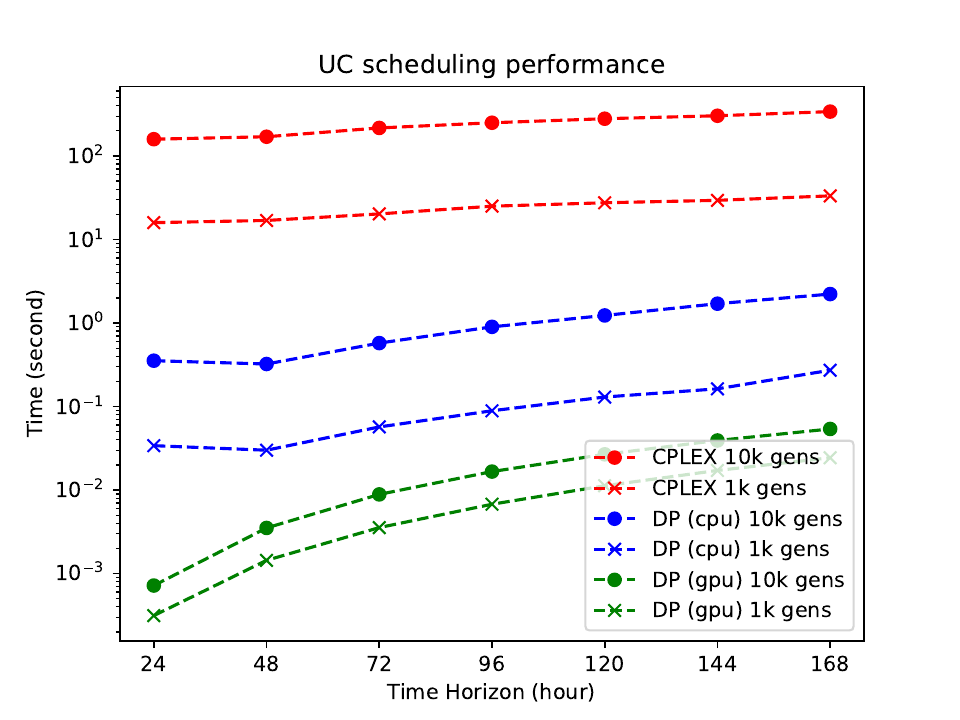}
    \caption{Computation time for solving UC with CPLEX vs. DP on CPU and GPU.}
    \label{fig:dp-performance}
\end{figure}

The numerical demonstration of the DP algorithm is important for both the CPU and GPU, because the UC subproblems are solved over the ADMM iterations.
Figure \ref{fig:dp-performance} shows the computation time for solving the UC subproblems with our DP algorithms on the CPU and GPU. We also solve these subproblems with CPLEX on the CPU, as a benchmark.
We consider UC with 1,000 and 10,000 generators by increasing the number of time periods from 24 to 168.
The figure shows that the DP algorithm solves the 1,000-generator subproblem about 100 times faster than CPLEX. Running on the GPU further accelerates the solution time more than 10,000 times faster.
While the solution time of the DP algorithm gradually increases with the number of time periods, we observe that the DP algorithm solves the 10,000-generator subproblem more than 10,000 times faster than CPLEX for 24 time periods.

\subsection{Performance of the Component-Based ADMM for UC-ACOPF}

In this subsection we demonstrate the numerical results from the component-based ADMM algorithm for solving UC-ACOPF with the test instances.
Table \ref{tab:admm} reports the penalty parameter values and numerical performance for each test instance. 
The primal infeasibility is computed as the maximum norm of the constraint violations of \eqref{eq:uc-acopf-ext}.

We observe that running the algorithm on the GPU is slower than that on the CPU for the small instances (case9 and case30), resulting in the speedup of 0.4 and 0.9, respectively.
However, significant speedups are achieved by running the algorithm on the GPU for the larger instances (case118 and case300).
For the case300 instance, the GPU accelerates the ADMM steps by 4.7 times.
In particular, the per-iteration CPU time significantly increases from 0.01 to 0.33 second with the increasing number of components (case9 vs. case300, respectively), whereas the per-iteration GPU time increases from 0.03 to 0.05 second for the same instances.
The main reason is that the GPU excels at performing a large number of simple computations simultaneously, making it ideal for tasks that require parallel processing.

\subsection{Optimality Gap from the Component-Based ADMM}

\begin{table}[!ht]
    \centering
    \begin{tabular}{l|rrrrr}
             &             &             & \multicolumn{2}{c}{Time (s)} \\
        Case & Lower Bound & Upper Bound & Root Note & Total \\
        \hline
        case9   & 8.6893e+4 & 9.0263e+4 & 4.0   & 61.1 \\
        case30  & 1.0433e+4 & 1.0442e+4 & 7.2   & 160.0 \\
        case118 & 2.1293e+6 & 2.1388e+6 & 98.1  & 1515.3 \\
        case300 & 1.1776e+7 & $\infty$  & 134.9 & 3600.0
    \end{tabular}
    \caption{Computational performance of solving the continuous relaxation by using \texttt{Juniper.jl} with \texttt{Ipopt}}
    \label{tab:juniper}
\end{table}

The quality of the solutions resulting from our ADMM algorithm is partly measured by the primal infeasibility reported in Table \ref{tab:admm}.
To assess the optimality gap of the objective values, we use \texttt{Juniper.jl} \cite{juniper}, a Julia-based branch-and-bound solver that can heuristically find the lower and upper bounds of MINLP problem.
Table \ref{tab:juniper} reports the lower and upper bounds of the test instances with the solution time. The solver is configured to use \texttt{Ipopt} and \texttt{HiGHS} for nonlinear optimization and mixed-integer programming solvers, respectively. We set the 1-hour time limit for each instance.
Not surprisingly, the solution time of the branch-and-bound solver significantly increases with the instance size for both the root node solution and total solution. No upper bound was found for the largest instance with 300 buses within the time limit.

\begin{table}[!ht]
    \centering
    \begin{tabular}{l|rr}
        Case & \texttt{Juniper.jl} & Our ADMM \\
        \hline
        case9 & 3.73 & 4.81 \\
        case30 & 0.08 & 1.60 \\
        case118 & 0.44 & 2.01 \\
        case300 & NA & 3.55
    \end{tabular}
    \caption{Optimality gap (\%) from \texttt{Juniper.jl} and the component-based ADMM.}
    \label{tab:gap}
\end{table}

We compute the optimality gap as $(\text{upper\_bound} - \text{lower\_bound})/\text{upper\_bound}$ and report the values in Table \ref{tab:gap}. Note that the objective values reported in Table \ref{tab:admm} are the upper bounds from the ADMM algorithm.
Table \ref{tab:gap} shows that the optimality gaps obtained by \texttt{Juniper.jl} are small in general for the instances solved. We observe that the gaps obtained by the ADMM are also small, but larger than those by \texttt{Juniper.jl}.

\section{Concluding Remarks and Future Work}
\label{sec:conclusion}

As shown and discussed in the previous work \cite{kim2021leveraging,kim2022accelerated}, the component-based ADMM algorithm \cite{mhanna2018adaptive} has been adapted to solve the UC-ACOPF problem GPU. 
To this end, we developed a novel dynamic programming approach to solving the UC subproblem with the time complexity of $O(T)$, resulting in the UC subproblem solution up to $10^5$ times faster than with CPLEX.
Our numerical experiments demonstrated the promising results from the component-based ADMM algorithm by utilizing the parallel processing of many small optimization subproblems on GPU. The results reported in this paper include the solution times on GPU and the optimality gaps obtained.

We have also identified the limitations from this work. A main limitation of the algorithm is that the numerical convergence is sensitive to the choice of the penalty parameter values (i.e., $\rho_\text{pq},\rho_\text{va},\rho_\text{uc}$). We also hypothesize that the parameter values also affect the solution quality and optimality gap.
Unlike branch-and-bound algorithms, the ADMM algorithm does not provide a lower bound or the optimality gap. GAIL - do you mea for the optimality gap?

Addressing such limitations poses important areas for future work.
For example, algorithmic approaches to updating the penalty parameters (e.g., adaptive rules \cite{sun2021two}, reinforcement learning \cite {zeng2022reinforcement}) need to be developed particularly for UC-ACOPF with the understanding of the connection to the solution quality.

\bibliographystyle{IEEEtran}
\bibliography{reference.bib}


\end{document}